\documentclass[11pt,english]{article}
\usepackage[T1]{fontenc}
\usepackage[latin9]{inputenc}
\usepackage{amsmath}
\usepackage{amssymb}

\makeatletter

\usepackage{color}

\usepackage{amsfonts}

\usepackage{a4wide}

\providecommand{\U}[1]{\protect\rule{.1in}{.1in}} \makeatletter
\newtheorem{theorem}{Theorem}

\newtheorem{corollary}[theorem]{Corollary}

\newtheorem{lemma}[theorem]{Lemma}

\newtheorem{proposition}[theorem]{Proposition}
\newtheorem{remark}{Remark}

\usepackage{babel}
\makeatother

\begin{document}

\title{Strongly-Representable Operators}

\author{M.D. Voisei%
\thanks{Towson University, Department of Mathematics, Towson, Maryland, USA,
email: \texttt{{mvoisei@towson.edu}.}%
}~ and C. Z\u{a}linescu%
\thanks{University {}``Al.I.Cuza'' Ia\c{s}i, Faculty of Mathematics, 700506-Ia\c{s}i,
Romania and Institute of Mathematics Octav Mayer, Ia\c{s}i, Romania,
email: \texttt{{zalinesc@uaic.ro}.}%
}}

\date{{}}

\maketitle
\begin{abstract}
Recently in \cite{alves 2008} a new class of maximal monotone operators
has been introduced. In this note we study domain range properties
as well as connections with other classes and calculus rules for these
operators we called strongly-representable. While not every maximal
monotone operator is strongly-representable, every maximal monotone
NI operator is strongly-representable, and every strongly representable
operator is locally maximal monotone, maximal monotone locally, and
ANA. As a consequence the conjugate of the Fitzpatrick function of
a maximal monotone operator is not necessarily a representative function. 
\end{abstract}

\section{Introduction}

Let $X$ be a non trivial Banach space and $X^{\ast}$ its topological
dual; set $Z:=X\times X^{\ast}$ which is a Banach space with respect
to the norm $\left\Vert (x,x^{\ast})\right\Vert :=\big(\left\Vert x\right\Vert ^{2}+\left\Vert x^{\ast}\right\Vert ^{2}\big)^{1/2}$.
We denote by $s$ the strong topology and by $Z^{*}:=X^{*}\times X^{**}$
the dual of $Z$.

For $z:=(x,x^{\ast})\in Z$ set $c(z):=\left\langle x,x^{\ast}\right\rangle $.
Consider \[
\mathcal{F}:=\mathcal{F}(Z):=\left\{ f\in\Lambda(Z)\mid f(z)\geq c(z)\ \forall z\in Z\right\} ,\quad\mathcal{F}_{s}:=\mathcal{F}_{s}(Z):=\mathcal{F}(Z)\cap\Gamma_{s}(Z),\]
 where for a locally convex space $(E,\tau)$, $\Lambda(E)$ denotes
the class of proper convex functions $f:E\rightarrow\overline{\mathbb{R}}$
and $\Gamma_{\tau}(E)$ is the class of those $f\in\Lambda(E)$ which
are $\tau$--lower semi\-continuous (lsc for short). The elements
of $\mathcal{F}(Z)$ are called representative functions in $Z$.
The classes $\mathcal{F}(Z^{\ast})$, $\mathcal{F}_{s}(Z^{\ast})$
are defined similarly.

It is known that when $f\in\mathcal{F}(Z)$ the set \[
M_{f}:=\left\{ z\in Z\mid f(z)\leq c(z)\right\} =\left\{ z\in Z\mid f(z)=c(z)\right\} =[f=c]\]
 is monotone, that is, $c(z-z^{\prime})\geq0$ for all $z,z^{\prime}\in M_{f}$.
For $z_{1}:=(x_{1},x_{1}^{\ast})$, $z_{2}:=(x_{2},x_{2}^{\ast})\in Z$
we set \[
\left\langle z_{1},z_{2}\right\rangle :=z_{1}\cdot z_{2}:=\left\langle x_{1},x_{2}^{\ast}\right\rangle +\left\langle x_{2},x_{1}^{\ast}\right\rangle .\]
 Note the following useful relations: \[
c(z_{1}+z_{2})=c(z_{1})+\left\langle z_{1},z_{2}\right\rangle +c(z_{2}),\quad c(z)=c(-z)=\tfrac{1}{2}\left\langle z,z\right\rangle \quad\forall z_{1},z_{2},z\in Z.\]

For $z=(x,x^{\ast})\in Z$, $\alpha>0$ and $g:Z\rightarrow\overline{\mathbb{R}}$
we denote by $g_{z}$ and $g_{\alpha}$ the functions defined on $Z$
by \[
g_{z}(w):=g(z+w)-c(z+w)+c(w),\quad g_{\alpha}(w):=\alpha g\left(y,\alpha^{-1}y^{\ast}\right),\]
 for $w:=(y,y^{\ast})\in Z;$ hence \begin{gather}
g_{z}(w)-c(w)=g(z+w)-c(z+w)\quad\forall z,w\in Z,\label{z0-a}\\
g_{\alpha}(w)-c(w)=\alpha\left[g(w_{\alpha})-c(w_{\alpha})\right]\quad\forall\alpha>0,~\forall w\in Z,\label{z0-b}\end{gather}
 where $w_{\alpha}:=(y,\alpha^{-1}y^{\ast})$ for $w=(y,y^{\ast})$.
It follows that \begin{gather*}
f\in\mathcal{F}(Z)\Rightarrow\left[f_{\alpha},\ f_{z}\in\mathcal{F}(Z)\quad\forall\alpha>0,\ \forall z\in Z\right],\\
f\in\mathcal{F}_{s}(Z)\Rightarrow\left[f_{\alpha},\
f_{z}\in\mathcal{F}_{s}(Z)\quad\forall\alpha>0,\ \forall z\in Z\right],\end{gather*}
 and\begin{equation}
M_{f_{z}}=M_{f}-z,\quad M_{f_{\alpha}}=\left\{ (x,\alpha x^{\ast})\mid(x,x^{\ast})\in M_{f}\right\} \label{z12}\end{equation}
 for every $f\in\mathcal{F}(Z)$, $z\in Z$ and $\alpha>0$. In the
sequel for $g:Z\rightarrow\overline{\mathbb{R}}$ a proper function
$g^{\ast}$ denotes its usual conjugate, while $\partial g$ is its
usual sub\-differential, that is, $g^{\ast}:Z^{\ast}=X^{\ast}\times X^{\ast\ast}\rightarrow\overline{\mathbb{R}}$
and $\partial g:Z\rightrightarrows Z^{\ast}$, the pairing between
$Z$ and $Z^{\ast}$ being given by\[
\left\langle (x,x^{\ast}),(u^{\ast},u^{\ast\ast})\right\rangle :=\left\langle x,u^{\ast}\right\rangle +\left\langle x^{\ast},u^{\ast\ast}\right\rangle \quad\forall(x,x^{\ast})\in X\times X^{\ast},\ (u^{\ast},u^{\ast\ast})\in X^{\ast}\times X^{\ast\ast}.\]

Let $\widehat{x}$ be the image $J(x)$ of $x\in X$, where $J$ is
the canonical injection of $X$ into $X^{\ast\ast}$, $J:X\rightarrow X^{\ast\ast}$
with $J(x)(x^{\ast}):=\left\langle x,x^{\ast}\right\rangle $ for
$x^{\ast}\in X^{\ast}$. In the sequel we shall use $\widehat{z}$
for $(x^{\ast},\widehat{x})\in Z^{\ast}$ when $z=(x,x^{\ast})\in Z$.
Moreover, for $g:Z\rightarrow\overline{\mathbb{R}}$ we consider $g^{\square}:Z\rightarrow\overline{\mathbb{R}}$
defined by $g^{\square}(z):=g^{\ast}(\widehat{z});$ hence $g^{\square}$
is convex and $s\times w^{\ast}$--lsc.

For $M\subset X\times X^{\ast}$, its Fitzpatrick function $\varphi_{M}$
is defined as\[
\varphi_{M}(z)=\sup\{\langle z,w\rangle-c_{M}(w)\mid w\in Z\},\]
 where $c_{M}(z):=c(z)$ for $z\in M$ and $c_{M}(z):=\infty$ for
$z\in Z\setminus M$; in simpler words $\varphi_{M}(x,x^{\ast})=c_{M}^{\ast}(x^{\ast},\widehat{x})=c_{M}^{\square}(x,x^{\ast})$,
or $\varphi_{M}(z)=c_{M}^{\ast}(\widehat{z})=c_{M}^{\square}(z)$
for $z=(x,x^{\ast})\in Z$.

Let $g:X\times X^{\ast}\rightarrow\overline{\mathbb{R}}$ be a proper
function and $z:=(x,x^{\ast})\in Z$. Then \[
(g_{(x,x^{\ast})})^{*}(u^{\ast},u^{\ast\ast})=g^{\ast}\left(u^{\ast}+x^{\ast},u^{\ast\ast}+\widehat{x}\right)-\left\langle x,u^{\ast}\right\rangle -\left\langle x^{\ast},u^{\ast\ast}\right\rangle -\left\langle x,x^{\ast}\right\rangle \quad\forall(u^{\ast},u^{\ast\ast})\in X^{\ast}\times X^{\ast\ast},\]
 that is, \[
(g_{z})^{*}(w^{\ast})=g^{\ast}\left(w^{\ast}+\widehat{z}\right)-c(w^{\ast}+\widehat{z})+c(w^{\ast})\quad\forall w^{\ast}\in Z^{\ast},\]
 or equivalently\[
(g_{z})^{*}=(g^{*})_{\widehat{z}},\]
 and \[
\partial g_{z}(w)=\left\{ w^{\ast}\in Z^{\ast}\mid w^{\ast}+\widehat{z}\in\partial g(w+z)\right\} =\partial g(w+z)-\widehat{z}.\]
 In particular, $\operatorname{Im}\partial g_{z}=\operatorname{Im}\partial g-\widehat{z}$.
Moreover, for $\alpha>0$ we have \begin{gather*}
g_{\alpha}^{\ast}(u^{\ast},u^{\ast\ast})=\alpha g^{\ast}(\alpha^{-1}u^{\ast},u^{\ast\ast}),\\
(u^{\ast},u^{\ast\ast})\in\partial g_{\alpha}(x,x^{\ast})\Leftrightarrow(\alpha^{-1}u^{\ast},u^{\ast\ast})\in\partial g(x,\alpha^{-1}x^{\ast}).\end{gather*}

Let us consider the more restrictive classes \[
\mathcal{G}:=\mathcal{G}(Z):=\left\{ f\in\mathcal{F}(Z)\mid f^{\ast}(z^{\ast})\geq c(z^{\ast})~\forall z^{\ast}\in Z^{\ast}\right\} ,\quad\
\mathcal{G}_{s}:=\mathcal{G}_{s}(Z):=\mathcal{G}(Z)\cap\Gamma_{s}(Z).\]
 The classes $\mathcal{G}(Z^{\ast})$, $\mathcal{G}_{s}(Z^{\ast})$
are defined similarly.

Using the formulas above for $g_{z}^{\ast}$ and $g_{\alpha}^{\ast}$
we get\begin{equation}
f\in\mathcal{G}_{s}(Z)\Rightarrow\left[f_{\alpha},\
f_{z}\in\mathcal{G}_{s}(Z)\ \ \forall\alpha>0,\ \forall z\in Z\right].\label{z13}\end{equation}

An operator $M$ is called \emph{strongly-representable} in $Z$ whenever
there is $f\in\mathcal{G}_{s}(Z)$ such that $M=M_{f}$. In this case
$f$ is called a strong-representative of $M$.

It has been proven in \cite[Theorem 4.2]{alves 2008} that every strongly-representable
operator is maximal monotone. In this paper we show that not every
maximal monotone operator is strongly-representable by providing the
property of convexity for the closure of the range; property that
distinguishes between these two classes.

Consider\[
h:X\times X^{\ast}\rightarrow\mathbb{R},\quad h(x,x^{\ast})=\tfrac{1}{2}\left\Vert (x,x^{\ast})\right\Vert ^{2}=\tfrac{1}{2}\left\Vert x\right\Vert ^{2}+\tfrac{1}{2}\left\Vert x^{\ast}\right\Vert ^{2}.\]
 Since the dual norm on $X^{\ast}\times X^{\ast\ast}$ is given by
$\left\Vert (u^{\ast},u^{\ast\ast})\right\Vert =\big(\left\Vert u^{\ast}\right\Vert ^{2}+\left\Vert u^{\ast\ast}\right\Vert ^{2}\big)^{1/2}$
we know that $h^{\ast}(u^{\ast},u^{\ast\ast})=\tfrac{1}{2}\left\Vert (u^{\ast},u^{\ast\ast})\right\Vert ^{2}$.
Notice that \begin{equation}
h\geq\pm c,\quad h^{\ast}\geq\pm c\label{z17}\end{equation}
 (on the respective spaces). Moreover,\[
\partial h(x,x^{\ast})=F_{X}(x)\times F_{X^{\ast}}(x^{\ast}),\]
 where $F_{X}:X\rightrightarrows X^{\ast}$ is the duality mapping
of $X$, that is, \[
F_{X}(x):=\partial\big(\tfrac{1}{2}\left\Vert \cdot\right\Vert ^{2}\big)(x)=\big\{x^{\ast}\in X^{\ast}\mid\left\Vert x\right\Vert ^{2}=\left\Vert x^{\ast}\right\Vert ^{2}=\left\langle x,x^{\ast}\right\rangle \big\}\quad\forall x\in X,\]
 and similarly for $F_{X^{\ast}}$. Note that\begin{equation}
\left\vert \left\langle z,z^{\prime}\right\rangle \right\vert \leq\left\Vert z\right\Vert \cdot\left\Vert z^{\prime}\right\Vert ,\quad\left\vert c(z)-c(z^{\prime})\right\vert \leq\tfrac{1}{2}\left\Vert z-z^{\prime}\right\Vert ^{2}+\left\Vert z^{\prime}\right\Vert \cdot\left\Vert z-z^{\prime}\right\Vert \quad\forall z,z^{\prime}\in Z.\label{z18}\end{equation}
 Taking $z^{\prime}=z$ in the first inequality or $z^{\prime}=0$
in the second we get $\left\vert c(z)\right\vert \leq\tfrac{1}{2}\left\Vert z\right\Vert ^{2}$
for $z\in Z$.

In the sequel a multi\-function $S:E\rightrightarrows F$ is identified
with its graph $\operatorname*{gph}S:=\{(x,y)\mid y\in S(x)\}$ (when
there is no risk of confusion); so $\operatorname*{dom}S:=\Pr_{E}(S)$
and $\operatorname{Im}S:=\Pr_{F}(S)$. Moreover, $S^{-1}:F\rightrightarrows E$
has $\operatorname*{gph}S^{-1}:=\{(y,x)\mid(x,y)\in\operatorname*{gph}S\}$.

When $E,F$ are (real) linear spaces, $A,B\subset E$ and $\alpha\in\mathbb{R}$
we set $A+B:=\{a+b\mid a\in A,\ b\in B\}$ and $\alpha A:=\{\alpha a\mid a\in A\}$
with $A+\emptyset:=\emptyset$ and $\alpha\emptyset:=\emptyset$.
For $S,T:E\rightrightarrows F$ and $\alpha\in\mathbb{R}$ the multi\-functions
$S+T:E\rightrightarrows F$ and $\alpha S:E\rightrightarrows F$ have
their graphs $\operatorname*{gph}(S+T):=\left\{ (x,y+v)\mid(x,y)\in\operatorname*{gph}S,\
(x,v)\in\operatorname*{gph}T\right\} $, that is, $(S+T)(x)=S(x)+T(x)$, and $\operatorname*{gph}(\alpha S):=\left\{ (x,\alpha y)\mid(x,y)\in\operatorname*{gph}S\right\} $,
that is, $(\alpha S)(x)=\alpha S(x)$. Hence $\operatorname*{dom}(S+T)=\operatorname*{dom}S\cap\operatorname*{dom}T$,
$\operatorname{Im}(S+T)\subset\operatorname{Im}S+\operatorname{Im}T$,
$\operatorname*{dom}(\alpha S)=\operatorname*{dom}S$, $\operatorname{Im}(\alpha S)=\alpha\operatorname{Im}S$;
therefore $\operatorname*{gph}(S+T)$ is (generally) different of
$\operatorname*{gph}S+\operatorname*{gph}T$ and $\operatorname*{gph}(\alpha S)$
is different of $\alpha\operatorname*{gph}S$.

As usual, for the subset $A$ of the normed vector space $X$ and
$x\in X$ we set $d(x,A):=\inf\left\{ \left\Vert x-u\right\Vert \mid u\in A\right\} $
with the convention $\inf\emptyset:=\infty:=+\infty$. 

\section{Domain-range properties}

\begin{proposition} \label{p-z1}Assume that $f\in\mathcal{F}(Z)$,
and $z_{1},z_{2}\in Z$ and $\varepsilon_{1},\varepsilon_{2}\geq0$
are such that $f(z_{1})\leq c(z_{1})+\varepsilon_{1}$ and $f(z_{2})\leq c(z_{2})+\varepsilon_{2}$.
Then \[
c(z_{1}-z_{2})=\left\langle x_{1}-x_{2},x_{1}^{\ast}-x_{2}^{\ast}\right\rangle \geq-2(\varepsilon_{1}+\varepsilon_{2}).\]

\end{proposition}

Proof. Indeed,\[
c\left(\tfrac{1}{2}z_{1}+\tfrac{1}{2}z_{2}\right)\leq f\left(\tfrac{1}{2}z_{1}+\tfrac{1}{2}z_{2}\right)\leq\tfrac{1}{2}f(z_{1})+\tfrac{1}{2}f(z_{2})\leq\tfrac{1}{2}\left(c(z_{1})+\varepsilon_{1}\right)+\tfrac{1}{2}\left(c(z_{2})+\varepsilon_{2}\right),\]
 whence $-\tfrac{1}{2}(\varepsilon_{1}+\varepsilon_{2})\leq\tfrac{1}{4}c(z_{1}-z_{2})$.
The conclusion follows.\hfill{}$\square$

\begin{proposition} \label{t-v}Let $f\in\mathcal{G}(Z)$. Then:

(i) For every $z\in Z$ one has\[
\inf_{w\in Z}\left(f_{z}(w)+h(w)\right)=-\min_{w^{\ast}\in Z^{\ast}}\left[\left(f^{\ast}(\widehat{z}+w^{\ast})-c(\widehat{z}+w^{\ast})\right)+\left(h^{\ast}(w^{\ast})+c(w^{\ast})\right)\right]=0.\]

(ii) For every $z\in Z$ there is $z^{\ast}\in M_{f^{\ast}}$ such
that $\widehat{z}-z^{\ast}\in\operatorname*{gph}(-F_{X^{\ast}})$
and $\Vert\widehat{z}-z^{\ast}\Vert^{2}\leq2(f^{\ast}(\widehat{z})-c(\widehat{z}))$.
Moreover \begin{equation}
(\sqrt{2}-1)\Vert\widehat{z}-z^{\ast}\Vert\leq d(\widehat{z},M_{f^{\ast}})\leq\sqrt{2\left(f^{\ast}(\widehat{z})-c(\widehat{z})\right)}=\sqrt{2\left(f^{\square}(z)-c(z)\right)}.\label{rpz1}\end{equation}

(iii) For every $\alpha>0$, $\operatorname{Im}\big(\left(M_{f^{\ast}}\right)^{-1}+\alpha(F_{X^{\ast}})^{-1}\big)=X^{\ast}$.
\end{proposition}

Proof. (i) Taking into account the formulas related to $f_{z}$ we
may (and we do) assume that $z=0$. Because $f\geq c$ and $f^{\ast}\geq c$,
from (\ref{z17}) we obtain that $f+h\geq0$ and $f^{\ast}+h^{\ast}\geq0$.
Since $f$ is convex and $h$ is finite, convex and continuous on
$Z$, using the Fenchel duality theorem (see f.i. \cite[Cor. 2.8.5]{Zalinescu:02})
we obtain that \[
0\leq\inf(f+h)=-\min_{z^{\ast}\in Z^{\ast}}\left[f^{\ast}(z^{\ast})+h^{\ast}(-z^{\ast})\right]=-\min_{z^{\ast}\in Z^{\ast}}\left[f^{\ast}(z^{\ast})+h^{\ast}(z^{\ast})\right]=-\inf(f^{\ast}+h^{\ast})\leq0.\]
 The conclusion of (i) follows, because $f_{z}\in\mathcal{G}(Z)$
whenever $f\in\mathcal{G}(Z)$.

(ii) Fix $z\in Z$. From (i) we get $z^{\ast}\in Z^{\ast}$ such that
$\left[f^{\ast}(z^{\ast})-c(z^{\ast})\right]+\left[h^{\ast}(z^{\ast}-\widehat{z})+c(z^{\ast}-\widehat{z})\right]=0$.
Because the terms in square brackets are non negative, we obtain that
$f^{\ast}(z^{\ast})-c(z^{\ast})=0$, that is, $z^{\ast}\in M_{f^{\ast}}$,
and $h^{\ast}(z^{\ast}-\widehat{z})+c(z^{\ast}-\widehat{z})=0$, that
is, $\widehat{z}-z^{\ast}\in\operatorname*{gph}(-F_{X^{\ast}})$.
Since $f^{\ast}(z^{\ast})=c(z^{\ast})$ we have $f^{\ast}(\widehat{z})\geq\varphi_{M_{f^{\ast}}}(\widehat{z})\geq\left\langle \widehat{z},z^{\ast}\right\rangle -c(z^{\ast})$
(for more details see \cite[Remark 3.6]{Voisei-tscr}). Therefore\[
f^{\ast}(\widehat{z})-c(\widehat{z})\geq\left\langle \widehat{z},z^{\ast}\right\rangle -c(z^{\ast})-c(\widehat{z})=-c(\widehat{z}-z^{\ast})=h^{\ast}(z^{\ast}-\widehat{z})=\tfrac{1}{2}\Vert\widehat{z}-z^{\ast}\Vert^{2}.\]
 Since $\delta:=d(\widehat{z},M_{f^{\ast}})\leq\left\Vert \widehat{z}-z^{\ast}\right\Vert $,
the second inequality in relation (\ref{rpz1}) holds. Since $M_{f^{\ast}}$
is monotone we have that\begin{align*}
0 & \leq c\left(z^{\ast}-w^{\ast}\right)=c(z^{\ast}-\widehat{z})+\left\langle z^{\ast}-\widehat{z},\widehat{z}-w^{\ast}\right\rangle +c(\widehat{z}-w^{\ast})\\
 & \leq-\tfrac{1}{2}\Vert\widehat{z}-z^{\ast}\Vert^{2}+\left\Vert z^{\ast}-\widehat{z}\right\Vert \cdot\left\Vert \widehat{z}-w^{\ast}\right\Vert +\tfrac{1}{2}\left\Vert \widehat{z}-w^{\ast}\right\Vert ^{2}\end{align*}
 for every $w^{\ast}\in M_{f^{\ast}}$. It follows that $0\leq-\Vert\widehat{z}-z^{\ast}\Vert^{2}+2\delta\left\Vert z^{\ast}-\widehat{z}\right\Vert +\delta^{2}$,
whence $\Vert\widehat{z}-z^{\ast}\Vert\leq(1+\sqrt{2})\delta$. Therefore,
the first inequality in (\ref{rpz1}) holds, too.

(iii) Replacing, if necessary, $f$ by $f_{\alpha}$, we may assume
that $\alpha=1$. Let $u^{\ast}\in X^{\ast}$. Applying (ii) for $z=(0,u^{\ast})$
we get $z^{\ast}=(x^{\ast},x^{\ast\ast})\in M_{f^{\ast}}$ such that
$u^{\ast}-x^{\ast}\in(F_{X^{\ast}})^{-1}(x^{\ast\ast})$. The conclusion
follows. \hfill{}$\square$

\begin{remark} \label{rem1}From assertion (i) of the preceding proposition
we have that $f\in\mathcal{G}_{s}(Z)$ implies $f\in\mathcal{F}_{s}(Z)$
and $\inf(f_{z}+h)=0$ for every $z\in Z$. \end{remark}

\begin{remark} \label{rem2}The first part of assertion (ii) of the
previous proposition can be interpreted as\[
\widehat{Z}:=X^{\ast}\times J(X)\subset\operatorname*{gph}M_{f^{\ast}}+\operatorname*{gph}(-F_{X^{\ast}}),\]
 and is a generalization to non-reflexive spaces for the {}``$-J$\textquotedblright\ criterion
for the maximality of operators in reflexive spaces (see \cite{Simons:98}).
In reflexive spaces, an operator is maximal monotone iff it is strongly-representable,
a situation that is no longer valid in the non-reflexive context in
the sense that there exist maximal monotone operators that are not
strongly-representable as we will see in the sequel. The second part
of assertion (ii) extends \cite[Lem. 2.3]{Penot/Zalinescu:06} to
the non-reflexive case. \end{remark}

A partial converse of Proposition \ref{t-v} follows.

\begin{proposition} \label{indep-cond} If $f:Z\rightarrow\overline{\mathbb{R}}$
is such that $\inf\nolimits _{w\in Z}\left(f_{z}(w)+h(w)\right)=0$
for every $z\in Z$ then $f\geq c;$ if moreover $f$ is convex then
$f\in\mathcal{F}(Z)$ and $f^{\ast}(z^{\ast})\geq c(z^{\ast})$ for
every $z^{\ast}\in\widehat{Z}+\operatorname*{gph}(-F_{X^{\ast}})$.
\end{proposition}

Proof. The condition $\inf\left(f_{z}+h\right)=0$ for every $z\in Z$
implies\[
f_{z}(w)+h(w)=f(z+w)-c(z+w)+h(w)+c(w)\geq0\quad\forall z,w\in Z.\]
 Taking $w=0$ we get $f\geq c$ in $Z$.

Assume now that $f$ is convex; then necessarily $f\in\Lambda(Z)$,
and so $f\in\mathcal{F}(Z)$. Again, the fundamental duality formula
yields\[
\inf_{w\in Z}\left(f_{z}(w)+h(w)\right)=-\min_{z^{\ast}\in Z^{\ast}}\left[\left(f^{\ast}(\widehat{z}+z^{\ast})-c(\widehat{z}+z^{\ast})\right)+\left(h^{\ast}(z^{\ast})+c(z^{\ast})\right)\right]=0,\]
 which implies $f^{\ast}(z^{\ast})\geq c(z^{\ast})$ for every $z^{\ast}\in\widehat{Z}+\operatorname*{gph}(-F_{X^{\ast}})$,
since $[h^{\ast}+c=0]=\operatorname*{gph}(-F_{X^{\ast}})$. \hfill{}$\square$

\begin{theorem} \label{t-z5}Let $f\in\Gamma_{s}(Z)$ be such that
$\inf\nolimits _{w\in Z}\left(f_{z}(w)+h(w)\right)=0$ for every $z\in Z$.
Then $M_{f}$ is monotone and nonempty, and \begin{equation}
d\left((x,x^{\ast}),M_{f}\right)\leq2\sqrt{f(x,x^{\ast})-\left\langle x,x^{\ast}\right\rangle }\quad\forall(x,x^{\ast})\in X\times X^{\ast}.\label{z11}\end{equation}

\end{theorem}

Proof. From Proposition \ref{indep-cond} we have that $f\in\mathcal{F}(Z)$,
and so $M_{f}$ is monotone. Fix $z:=(x,x^{\ast})\in X\times X^{\ast}$.
If $f(z)=\infty$ or $f(z)=c(z)$ it is nothing to prove. So let $\varepsilon:=f(z)-c(z)\in(0,\infty)$
and set $\varepsilon_{0}:=\varepsilon$, $z_{0}:=z$. Fix $\beta\in(1,\infty)$
and $\gamma\in(2,\infty)$ and consider a sequence $(\varepsilon_{n})_{n\geq0}\subset(0,\infty)$
satisfying \begin{equation}
4\varepsilon_{n}+6\varepsilon_{n+1}\leq\gamma^{2}\varepsilon_{n}\quad\forall n\geq0\qquad\text{and}\qquad\sum_{n\geq0}\sqrt{\varepsilon_{n}}<\beta\sqrt{\varepsilon}.\label{z10}\end{equation}
 Because $\inf(f_{z_{0}}+h)=0$, there exists $z_{1}\in Z$ such that
\[
f_{z_{0}}\left(z_{1}-z_{0}\right)+h\left(z_{1}-z_{0}\right)\leq\varepsilon_{1}.\]
 Using the definition of $f_{z}$ given in (\ref{z0-a}), we get \begin{gather}
0\leq f(z_{1})-c(z_{1})=f_{z_{0}}\left(z_{1}-z_{0}\right)-c(z_{1}-z_{0})\leq\varepsilon_{1},\nonumber \\
0\leq\tfrac{1}{2}\left\Vert z_{1}-z_{0}\right\Vert ^{2}+c(z_{1}-z_{0})\leq\varepsilon_{1}.\label{z9-b}\end{gather}
 Using Proposition \ref{p-z1} we obtain that $c(z_{1}-z_{0})\geq-2(\varepsilon_{0}+\varepsilon_{1})$,
and so, by (\ref{z9-b}), \[
\left\Vert z_{1}-z_{0}\right\Vert ^{2}\leq2\varepsilon_{1}+4(\varepsilon_{0}+\varepsilon_{1})=4\varepsilon_{0}+6\varepsilon_{1}\leq\gamma^{2}\varepsilon_{0},\]
 whence \[
\left\Vert z_{1}-z_{0}\right\Vert \leq\gamma\sqrt{\varepsilon_{0}}.\]
 Continuing this procedure we obtain a sequence $\left(z_{n}\right)_{n\geq0}\subset Z$
such that \[
f(z_{n})\leq c(z_{n})+\varepsilon_{n},\quad\left\Vert z_{n+1}-z_{n}\right\Vert \leq\gamma\sqrt{\varepsilon_{n}}\quad\forall n\geq0.\]
 From (\ref{z10}) we obtain that \[
\sum_{n\geq0}\left\Vert z_{n+1}-z_{n}\right\Vert \leq\gamma\sum_{n\geq0}\sqrt{\varepsilon_{n}}<\gamma\beta\sqrt{\varepsilon}.\]
 It follows that the sequence $\left(z_{n}\right)_{n\geq0}$ is strongly
convergent to some $z_{\varepsilon}\in Z$ and $\left\Vert z-z_{\varepsilon}\right\Vert \leq\gamma\beta\sqrt{\varepsilon}$.
Since $f$ is $s$--lsc and $\varepsilon_{n}\rightarrow0$, from the
inequality $f(z_{n})\leq c(z_{n})+\varepsilon_{n}$ we obtain \[
c(z_{\varepsilon})\leq f(z_{\varepsilon})\leq\liminf f(z_{n})\leq\lim(c(z_{n})+\varepsilon_{n})=c(z_{\varepsilon}).\]
 Therefore, $f(z_{\varepsilon})=c(z_{\varepsilon})$, that is, $z_{\varepsilon}\in M_{f}$.
Moreover, $d\left(z,M_{f}\right)\leq\gamma\beta\sqrt{\varepsilon}$.
Since $\beta>1$ and $\gamma>2$ are arbitrary we have that $d\left(z,M_{f}\right)\leq2\sqrt{\varepsilon}$,
that is, (\ref{z11}) holds. \hfill{}$\square$

\medskip{}

As a consequence of the previous theorem, every strongly-representable
operator has the following Br{\o}ndsted--Rockafellar property.
For other results of this type see \cite{alves 2008}.

\begin{corollary} \label{BR}Let $f\in\Gamma_{s}(Z)$ be such that
$\inf\nolimits _{w\in Z}\left(f_{z}(w)+h(w)\right)=0$, for every
$z\in Z$. For every $\varepsilon>0$ and every $(x,x^{\ast})\in X\times X^{\ast}$
with $f(x,x^{\ast})<\left\langle x,x^{\ast}\right\rangle +\varepsilon$
there exists $(x_{\varepsilon},x_{\varepsilon}^{\ast})\in M_{f}$
such that $\Vert x-x_{\varepsilon}\Vert^{2}+\Vert x^{\ast}-x_{\varepsilon}^{\ast}\Vert^{2}<4\varepsilon$.
\end{corollary}

The next result corresponds to \cite[Prop. 2]{Penot/Zalinescu:08}
(established in reflexive Banach spaces).

\begin{corollary} \label{c-z1}Let $f\in\mathcal{G}_{s}(Z)$ and
$\gamma>4$. Then for every $(x,x^{\ast})\in X\times X^{\ast}$ and
every $\alpha>0$ there exists $(x_{\alpha},x_{\alpha}^{\ast})\in M_{f}$
such that \begin{equation}
\left\Vert x_{\alpha}-x\right\Vert ^{2}+\alpha^{2}\left\Vert x_{\alpha}^{\ast}-x^{\ast}\right\Vert ^{2}\leq\gamma\alpha\left(f(x,x^{\ast})-\left\langle x,x^{\ast}\right\rangle \right).\label{z14}\end{equation}

\end{corollary}

Proof. If $(x,x^{\ast})\notin\operatorname*{dom}f$ we can take arbitrary
$(x_{\alpha},x_{\alpha}^{\ast})\in M_{f}$, while if $f(x,x^{\ast})=\left\langle x,x^{\ast}\right\rangle $
we take $(x_{\alpha},x_{\alpha}^{\ast})=(x,x^{\ast})$ for every $\alpha>0$.
So let $(x,x^{\ast})$ be such that $0<f(x,x^{\ast})-\left\langle x,x^{\ast}\right\rangle <\infty$
and fix $\alpha>0$. By (\ref{z13}) we have that $f_{\alpha}\in\mathcal{G}_{s}^{{}}(Z)$;
moreover, \[
f_{\alpha}(x,\alpha x^{\ast})-\left\langle x,\alpha x^{\ast}\right\rangle =\alpha\left(f(x,x^{\ast})-\left\langle x,x^{\ast}\right\rangle \right)\in(0,\infty).\]
 Applying Theorem \ref{t-z5} for $f_{\alpha}$ and $(x,\alpha x^{\ast})$
we get $(x_{\alpha},x_{\alpha}^{\ast})\in M_{f}$ (that is, $(x_{\alpha},\alpha x_{\alpha}^{\ast})\in M_{f_{\alpha}}$)
such that (\ref{z14}) holds. \hfill{}$\square$

\begin{corollary} \label{c-z2}Let $f\in\mathcal{G}_{s}(Z)$. Then
\[
\operatorname*{cl}(\operatorname*{dom}M_{f})=\operatorname*{cl}\left(\Pr\nolimits _{X}(\operatorname*{dom}f)\right),\qquad\operatorname*{cl}(\operatorname{Im}M_{f})=\operatorname*{cl}\left(\Pr\nolimits _{X^{\ast}}(\operatorname*{dom}f)\right).\]
 In particular $\operatorname*{cl}(\operatorname*{dom}M_{f})$ and
$\operatorname*{cl}(\operatorname{Im}M_{f})$ are convex sets. Here
{}``$\operatorname*{cl}$\textquotedblright\ stands for the closure
with respect to the strong topology. \end{corollary}

Proof. The inclusions $\operatorname*{dom}M_{f}\subset\Pr_{X}(\operatorname*{dom}f)$
and $\operatorname{Im}M_{f}\subset\Pr_{X^{\ast}}(\operatorname*{dom}f)$
are obvious. Let $x^{\ast}\in\Pr_{X^{\ast}}(\operatorname*{dom}f)$;
then $(x,x^{\ast})\in\operatorname*{dom}f$ for some $x\in X$. Applying
Corollary \ref{c-z1} (with $\gamma>4$), for $\alpha>0$ we get $(x_{\alpha},x_{\alpha}^{\ast})\in M_{f}$
satisfying (\ref{z14}). Therefore, $x_{\alpha}^{\ast}\in\operatorname{Im}M_{f}$
and $\alpha\left\Vert x_{\alpha}^{\ast}-x^{\ast}\right\Vert ^{2}\leq\gamma\left(f(x,x^{\ast})-\left\langle x,x^{\ast}\right\rangle \right)$.
Hence $\lim_{\alpha\rightarrow\infty}x_{\alpha}^{\ast}=x^{\ast}$,
which shows that $x^{\ast}\in\operatorname*{cl}(\operatorname{Im}M_{f})$.

If $x\in\Pr_{X}(\operatorname*{dom}f)$, $(x,x^{\ast})\in\operatorname*{dom}f$
for some $x^{\ast}\in X^{\ast}$. Taking for $\alpha>0$ $(x_{\alpha},x_{\alpha}^{\ast})\in M_{f}$
satisfying (\ref{z14}), we have that $x_{\alpha}\in\operatorname*{dom}M_{f}$
and $\left\Vert x_{\alpha}-x\right\Vert ^{2}\leq\gamma\alpha\left(f(x,x^{\ast})-\left\langle x,x^{\ast}\right\rangle \right)$.
Hence $\lim_{\alpha\rightarrow0}x_{\alpha}=x$, which proves that
$x\in\operatorname*{cl}\left(\operatorname*{dom}M_{f}\right)$. \hfill{}$\square$

\begin{remark} \label{rem3}The previous result shows that, for a
strongly-representable operator, the strong closures of its domain
and range are convex. Since, in general, the closure of the range
for a maximal monotone operator is not necessarily convex (see e.g.
\cite{Fitz-phelps - 95}), this shows that not every maximal monotone
operator is strongly-representable. \end{remark}

\begin{remark} \label{rem4}Let $M$ be a maximal monotone operator
that is not strongly-representable. Then $M=[\varphi_{M}=c]$, $\varphi_{M}\in\mathcal{F}_{s}(Z)$
and if we assumed that $\varphi_{M}^{\ast}\geq c$ in $Z^{\ast}$
then $\varphi_{M}\in\mathcal{G}_{s}(Z)$ and $M$ would be strongly-representable;
a contradiction. Hence the inequality $\varphi_{M}^{\ast}\geq c$
fails in $Z^{\ast}$, that is, the conjugate of the Fitzpatrick function
of a maximal monotone operator is not necessarily a representative
function. \end{remark}

The next result has been proved in \cite[Theorem 4.2]{alves 2008}
for $f\in\mathcal{G}_{s}(Z)$. For convenience we provide the reader
with a short proof.

\begin{theorem} \label{SR-MM} Let $f\in\Gamma_{s}(Z)$ be such that
$\inf\nolimits _{w\in Z}\left(f_{z}(w)+h(w)\right)=0$ for every $z\in Z$.
Then $M_{f}$ is maximal monotone in $Z$. In particular every strongly-representable
operator is maximal monotone. \end{theorem}

Proof. Let $z_{0}$ be monotonically related to $M_{f}$. Replacing
$f$ by $f_{z_{0}}$ if necessary, we may assume without loss of generality
that $z_{0}=0$, that is\begin{equation}
c(z)\geq0\quad\ \forall z\in M_{f}.\label{mrt}\end{equation}
 From $\inf(f+h)=0$, there is $z_{n}\in Z$ such that $f(z_{n})+h(z_{n})<1/n^{2}$,
for every $n\geq1$. The function $f+h$ is coercive. Indeed, fixing
some $\overline{z}^{\ast}\in\operatorname*{dom}f^{\ast}$ we have
that \[
f(z)+h(z)\geq\tfrac{1}{2}\left\Vert z\right\Vert ^{2}+\left\langle z,\overline{z}^{\ast}\right\rangle -f^{\ast}(\overline{z}^{\ast})\geq\tfrac{1}{2}\left\Vert z\right\Vert ^{2}-\left\Vert z\right\Vert \left\Vert \overline{z}^{\ast}\right\Vert -f^{\ast}(\overline{z}^{\ast})\quad\forall z\in Z.\]
 Therefore, the sequence $(z_{n})_{n\geq1}$ is bounded. Since $f\ge c$
and $h\ge-c$ we obtain that $f(z_{n})<c(z_{n})+1/n^{2}$ and $h(z_{n})+c(z_{n})\leq1/n^{2}$.
Applying Corollary \ref{BR} for $z_{n}$, $f$ and $\varepsilon=1/n^{2}$
we get $u_{n}\in M_{f}$ such that $\Vert u_{n}-z_{n}\Vert<2/n$ for
$n\geq1$.

According to (\ref{mrt}) and (\ref{z18}) we get\begin{align*}
\Vert z_{n}\Vert^{2} & =2h(z_{n})\leq-2c(z_{n})+2n^{-2}\leq-2c(u_{n})+2\left\vert c(u_{n})-c(z_{n})\right\vert +2n^{-2}\\
 & \leq\left\Vert u_{n}-z_{n}\right\Vert ^{2}+2\left\Vert z_{n}\right\Vert \cdot\left\Vert u_{n}-z_{n}\right\Vert +2n^{-2}\leq6n^{-2}+4n^{-1}\left\Vert z_{n}\right\Vert \end{align*}
 for $n\geq1$. Since $(z_{n})$ is bounded we have that $\left\Vert z_{n}\right\Vert \rightarrow0$.
Letting $n\rightarrow\infty$ in $c(z_{n})\leq f(z_{n})<c(z_{n})+1/n^{2}$
and taking into account that $f\in\Gamma_{s}(Z)$ we get $z_{0}=0\in M_{f}$.
\hfill{}$\square$

\begin{remark} \label{rem5}The subdifferential $\partial\varphi$
of a the function $\varphi\in\Gamma_{s}(X)$ with $X$ a Banach space
is strongly-representable thus maximal monotone; a strong-representative
for $\partial\varphi$ is given by $f(x,x^{\ast})=\varphi(x)+\varphi^{\ast}(x^{\ast})$
for $x\in X$, $x^{\ast\ast}\in X^{\ast\ast}$.\end{remark}

\begin{corollary} \label{CL} Let $f\in\Gamma_{s}(Z)$ be such that
$\inf\nolimits _{w\in Z}\left(f_{z}(w)+h(w)\right)=0$, for every
$z\in Z$. Then \begin{equation}
f\ge{\operatorname*{cl}}_{w\times w^{\ast}}\! f\ge\varphi_{M_{f}}\geq c,\ \text{ in }\
Z,\label{ineq-s}\end{equation}
 $M_{f}=M_{\operatorname*{cl}_{w\times w^{\ast}}\! f}\subset[f^{\square}=c]$,
and $\inf\nolimits _{w\in Z}\left((\operatorname*{cl}_{w\times w^{\ast}}\! f)_{z}(w)+h(w)\right)=0$,
for every $z\in Z$. Here $\operatorname*{cl}_{w\times w^{\ast}}\! f$
stands for the greatest convex $w\times w^{\ast}-$lsc function majorized
by $f$ in $Z$. \end{corollary}

Proof. According to Theorem \ref{SR-MM}, $M_{f}$ is maximal monotone.
By \cite[Th.\ 2.4]{Fi}, if $z\in M_{f}$ then $z\in\partial f(z)$.
This implies $f^{\square}(z)=c(z)$ for every $z\in M_{f}$, that
is, $M_{f}\subset[f^{\square}=c]$ and so $f^{\square}\leq c_{M_{f}}$.
Hence $f\ge\operatorname*{cl}_{w\times w^{\ast}}\! f=f^{\square\square}\geq\varphi_{M_{f}}\geq c$
in $Z$. Therefore $0\leq(\operatorname*{cl}_{w\times w^{\ast}}\! f)_{z}+h\leq f_{z}+h$
and so $\inf\nolimits _{w\in Z}\left((\operatorname*{cl}_{w\times w^{\ast}}\! f)_{z}(w)+h(w)\right)=0$,
for every $z\in Z$.

From $f\geq\operatorname*{cl}_{w\times w^{\ast}}\! f\geq c$ we get
$M_{f}\subset M_{\operatorname*{cl}_{w\times w^{\ast}}\! f}$. Because
$M_{f}$ is maximal and $M_{\operatorname*{cl}_{w\times w^{\ast}}\! f}$
is monotone the equality ensues. \hfill{}$\square$

\medskip{}

As a direct consequence of the previous corollary and Proposition
\ref{t-v}, the next result shows that the representative of a strongly-monotone
operator can be picked to be lsc with respect to the $w\times w^{\ast}$
topology on $Z$. Set $\mathcal{G}_{w\times w^{\ast}}(Z):=\mathcal{G}(Z)\cap\Gamma_{w\times w^{\ast}}(Z)=\mathcal{G}(Z)\cap\Gamma_{s\times w^{\ast}}(Z)$.

\begin{corollary} \label{c-v1}For every $f\in\mathcal{G}_{s}(Z)$,
$\operatorname*{cl}_{w\times w^{\ast}}\! f\in\mathcal{G}_{w\times w^{\ast}}(Z)$
and $M_{f}=M_{\operatorname*{cl}_{w\times w^{\ast}}\! f}=M_{f^{\square}}$.
In particular $\{M_{f}\mid f\in\mathcal{G}_{s}(Z)\}=\{M_{f}\mid f\in\mathcal{G}_{w\times w^{\ast}}(Z)\}$.
Moreover, if $f$ is a strong representative of $M\subset Z$ then
so are $\operatorname*{cl}_{w\times w^{\ast}}\! f$ and $\varphi_{M}$.
\end{corollary}

Proof. As previously seen in Corollary \ref{CL}, $M_{f}=M_{\operatorname*{cl}_{w\times w^{\ast}}\! f}\subset[f^{\square}=c]$
and from $f^{*}\ge c$ we know that $f^{\square}\ge c$ and $M_{f^{\square}}=[f^{\square}=c]$
is monotone. Since $M_{f}$ is maximal monotone the equality holds.
Moreover, from (\ref{ineq-s}) we get $\varphi_{M_{f}}^{*}\ge({\operatorname*{cl}}_{w\times w^{\ast}}\! f)^{*}\ge f^{*}\ge c$
which proves that $\operatorname*{cl}_{w\times w^{\ast}}\! f$ and
$\varphi_{M_{f}}$ are strong representatives of $M_{f}$. \hfill{}$\square$

\strut

\begin{corollary} \label{CLS} For every $f\in\mathcal{G}(Z)$, $\overline{f}:=\operatorname*{cl}_{s}f\in\mathcal{G}_{s}(Z)$
and $M_{\overline{f}}=M_{f^{\square}}$ is a maximal monotone extension
of $M_{f}$.\end{corollary}

Proof. Since $f\geq c$ we have that $f\ge\overline{f}\geq c$, $M_{f}\subset M_{\overline{f}}=M_{\overline{f}^{\square}}=M_{f^{\square}}$,
and $\overline{f}\in\mathcal{G}_{s}(X\times X^{\ast})$ because $f^{*}=\overline{f}^{*}\ge c$
and $\overline{f}^{\square}=f^{\square}$.\hfill{}$\square$

\strut

An immediate consequence of the preceding results is the following
characterization of strongly-representable operators.

\begin{theorem} \label{c-v2}Let $M\subset X\times X^{*}$ be monotone.
The following are equivalent

(i) $M$ is strongly representable,

(ii) $\varphi_{M}\in{\cal G}(X\times X^{*})$ and $M$ is representable,
that is, there is $f\in{\cal F}_{w\times w^{*}}(X\times X^{*})$ such
that $M=M_{f}$,

(iii) $M$ is maximal monotone and $\varphi_{M}^{*}\ge c$.

\end{theorem}

Proof. The implication (i)$\Rightarrow$(ii) follows from Corollary
\ref{c-v1} with $f=\varphi_{M}$.

For (ii)$\Rightarrow$(iii) it suffices to prove that $M$ is maximal
monotone. According to \cite[Theorem 3.4]{Voisei-tscr}, condition
$M=M_{f}$ for some $f\in{\cal F}_{w\times w^{*}}(X\times X^{*})$
together with $\varphi_{M}\ge c$ imply that $M$ is maximal monotone.

If (iii) holds then $M=M_{\varphi_{M}}$ and $\varphi_{M}\ge c$.
Therefore $\varphi_{M}$ is a strong-representative of $M$.

\hfill{}$\square$

\section{Calculus rules for strongly-representable operators}

We base our argument on the construction used in \cite{Penot/Zalinescu:05}
and note that several results of Section 3 in \cite{Penot/Zalinescu:05}
are valid without the reflexivity assumption.

For $X,Y$ locally convex spaces and $F:X\times Y\rightrightarrows X^{\ast}\times Y^{\ast}$
we define the multi\-function $G:=G(F):X\rightrightarrows X^{\ast}$
by \[
\operatorname*{gph}G:=\{(x,x^{\ast})\in X\times X^{\ast}\mid\exists y^{\ast}\in Y^{\ast}:\ (x,0,x^{\ast},y^{\ast})\in F\}.\]

As noticed in \cite{Penot/Zalinescu:05}, $G$ is monotone whenever
$F$ is monotone.

In general, for a locally convex space $E$, we denote by $\mathcal{M}(E)$
the class of monotone subsets of $E\times E^{*}$ and by $\mathfrak{M}(E)$
the class of maximal monotone subsets of $E\times E^{*}$. Moreover,
we denote by $\operatorname*{aff}A$ and $\overline{\operatorname*{aff}}A$
the affine hull and the closed affine hull of $A\subset E$, respectively.

First consider the following slight generalization of \cite[Lemma 3.1]{Penot/Zalinescu:05}.

\begin{lemma} \label{lem4} Let $X,Y$ be separated locally convex
spaces.

(i) If $F\in\mathcal{M}(X\times Y)$ and $Y_{0}\subset Y$ is a closed
linear subspace such that\begin{equation}
F(x,y)=F(x,y)+\{0\}\times Y_{0}^{\bot}\quad\forall(x,y)\in X\times Y,\label{r2}\end{equation}
 then $\Pr_{Y}(\operatorname*{dom}\varphi_{F})\subset y+Y_{0}$ for
every $y\in\Pr_{Y}(F)$.

(ii) If $F\in\mathfrak{M}(X\times Y)$, then $\Pr_{Y}(\operatorname*{dom}\varphi_{F})\subset\overline{\operatorname*{aff}}(\Pr_{Y}(F))$.
\end{lemma}

Proof. (i) Fix $y\in\Pr_{Y}(F)$ that is $(x,y,x^{\ast},y^{\ast})\in F$
for some $(x,x^{\ast},y^{\ast})\in X\times X^{*}\times Y^{*}$. By
(\ref{r2}), for every $v^{\ast}\in Y_{0}^{\bot}$ we have $(x,y,x^{\ast},y^{\ast}+v^{\ast})\in F$.

For every $\overline{y}\in\Pr_{Y}(\operatorname*{dom}\varphi_{F})$
there exist $(\overline{x},\overline{x}^{\ast},\overline{y}^{\ast})\in X\times X^{\ast}\times Y^{\ast}$,
$\gamma\in\mathbb{R}$ such that $\varphi_{F}(\overline{x},\overline{y},\overline{x}^{\ast},\overline{y}^{\ast})\le\gamma$.
From the definition of $\varphi_{F}$ we have\begin{align*}
\gamma & \geq\left\langle (\overline{x},\overline{y}),(x^{\ast},y^{\ast}+v^{\ast})\right\rangle +\left\langle (x,y),(\overline{x}^{\ast},\overline{y}^{\ast})\right\rangle -\left\langle (x,y),(x^{\ast},y^{\ast}+v^{\ast})\right\rangle \\
 & =\left\langle \overline{x}-x,x^{\ast}\right\rangle +\left\langle \overline{y},y^{\ast}\right\rangle +\left\langle x,\overline{x}^{\ast}\right\rangle +\left\langle y,\overline{y}^{\ast}-y^{\ast}\right\rangle +\left\langle \overline{y}-y,v^{\ast}\right\rangle ,\end{align*}
 which provides us with \[
\left\langle y-\overline{y},v^{\ast}\right\rangle \geq0\quad\forall v^{\ast}\in Y_{0}^{\bot}.\]
 This implies that $\overline{y}-y\in(Y_{0}^{\bot})^{\bot}=Y_{0}$.
Hence $\Pr_{Y}(\operatorname*{dom}\varphi_{F})\subset y+Y_{0}$.

(ii) Take $Y_{0}:=\overline{\operatorname*{aff}}(\Pr_{Y}(F))-y$ for
$y\in\Pr_{Y}(F)$ fixed. The operator $F+\Phi$ with $\operatorname*{gph}\Phi:=X\times(y+Y_{0})\times\left\{ 0\right\} \times Y_{0}^{\bot}$
is monotone and contains the maximal monotone operator $F$, so it
coincides with $F$, from which (\ref{r2}) follows. From (i) we get
the conclusion.\hfill{}$\square$

\medskip{}

As in \cite{Penot/Zalinescu:05}, we use the notation $\operatorname*{ri}A$
for the topological interior of $A$ with respect to $\overline{\operatorname*{aff}}A$,
and ${}^{ic}A$ for the relative algebraic interior of $A$ with respect
to $\overline{\operatorname*{aff}}A$; thus $\operatorname*{ri}A$,
${}^{ic}A$ are empty if $\operatorname*{aff}A$ is not closed and
one always has $\operatorname*{ri}A\subset{}^{ic}A$. In the sequel,
we use the facts that for $C$ convex with $^{ic}C$ nonempty, we
have $\operatorname*{aff}C=\operatorname*{aff}(^{ic}C)$ and, \begin{equation}
^{ic}C\subset A\subset C\Longrightarrow\lbrack\operatorname*{aff}C=\operatorname*{aff}A\text{ and }{}^{ic}C={}^{ic}A].\label{suppl}\end{equation}

\begin{theorem} \label{lem12}Let $X,Y$ be Banach spaces and $f\in\mathcal{G}_{s}(X\times Y\times X^{\ast}\times Y^{\ast})$.

(i) If $0\in{}^{ic}(\Pr_{Y}(\operatorname*{dom}f))$ and $g:X\times X^{\ast}\rightarrow\overline{\mathbb{R}}$
is given by \begin{equation}
g(x,x^{\ast}):=\inf\{f(x,0,x^{\ast},y^{\ast})\mid y^{\ast}\in Y^{\ast}\},\ (x,x^{\ast})\in X\times X^{\ast},\label{r9}\end{equation}
 then $g\in\mathcal{G}(X\times X^{\ast})$, \begin{equation}
g^{\ast}(u^{\ast},u^{\ast\ast})=\min\{f^{\ast}(u^{\ast},v^{\ast},u^{\ast\ast},0)\mid v^{\ast}\in Y^{\ast}\}\quad\forall\,(u^{\ast},u^{\ast\ast})\in X^{\ast}\times X^{\ast\ast},\label{r7bis}\end{equation}
 $\overline{g}=\operatorname*{cl}_{s}g\in\mathcal{G}_{s}(X\times X^{\ast})$
and \begin{equation}
G(M_{f})=M_{g}=M_{\overline{g}}=M_{g^{\square}}.\label{r10}\end{equation}
 Moreover $G(M_{f})$ is strongly representable and $\overline{g}$
is a strong representative of $G(M_{f})$; in particular $G(M_{f})$
is maximal monotone.

(ii) One has \begin{equation}
{}^{ic}({\Pr}_{Y}(\operatorname*{dom}f))={}^{ic}(\operatorname*{conv}({\Pr}_{Y}(M_{f})))={}^{ic}({\Pr}_{Y}(M_{f}))=\operatorname*{ri}({\Pr}_{Y}(M_{f}))={}^{ic}({\Pr}_{Y}(\operatorname*{dom}\varphi_{M_{f}})).\label{r14}\end{equation}

Therefore, if $0\in{}^{ic}(\Pr_{Y}(M_{f}))$ then $G(M_{f})$ is maximal
monotone. \end{theorem}

Proof. (i) First observe, from their definitions, that $g\geq c$
and $G(M_{f})\subset M_{g}$. To get (\ref{r10}) we follow the proof
of \cite[Lemma 3.2]{Penot/Zalinescu:05}; just observe that this time
the graph of $\mathcal{C}:X\times X^{\ast}\rightrightarrows X\times Y\times X^{\ast}\times X^{\ast}$
given by \[
\mathcal{C}(x,x^{\ast}):=\{x\}\times\{0\}\times\{x^{\ast}\}\times Y^{\ast},\quad(x,x^{\ast})\in X\times X^{\ast},\]
 is a closed linear subspace and $\mathcal{C}^{\ast}(x^{\ast},y^{\ast},x^{\ast\ast},y^{\ast\ast})=\{(x^{\ast},x^{\ast\ast})\}$
if $y^{\ast\ast}=0$, $\mathcal{C}^{\ast}(x^{\ast},y^{\ast},x^{\ast\ast},y^{\ast\ast})=\emptyset$
otherwise.

Notice that $g(x,x^{\ast})=\inf\{f(u,v,u^{\ast},v^{\ast})\mid(u,v,u^{\ast},v^{\ast})\in\mathcal{C}(x,x^{\ast})\}$
for $(x,x^{\ast})\in X\times X^{\ast}$ and\[
\operatorname*{dom}f-\operatorname*{Im}\mathcal{C}=X\times\Pr\nolimits _{Y}(\operatorname*{dom}f)\times X^{\ast}\times Y^{\ast},\]
 from which $0\in{}^{ic}(\operatorname*{dom}f-\operatorname*{Im}\mathcal{C})$.

By the fundamental duality formula (more precisely see \cite[Theorem 2.8.6 (v)]{Zalinescu:02})
we get (\ref{r7bis}). Since $f^{\ast}\geq c$, from (\ref{r7bis}),
we see that $g^{\ast}\geq c$, and so $g\in\mathcal{G}(X\times X^{\ast})$.

Since $g\in{\cal G}(X\times X^{*})$ we know by Corollary \ref{CLS}
that $M_{g}\subset M_{\overline{g}}=M_{g^{\square}}$ and $\overline{g}\in\mathcal{G}_{s}(X\times X^{\ast})$.
Therefore, according to Corollary \ref{c-v1} and again from (\ref{r7bis})
\[
M_{g^{\square}}=G(M_{f^{\square}})=G(M_{f})\subset M_{g}\subset M_{\overline{g}}=M_{g^{\square}}.\]
 Hence (\ref{r10}) holds.

(ii) Set $F:=M_{f}$. We first claim that \begin{equation}
{}^{ic}({\Pr}_{Y}(\operatorname*{dom}f))\subset{\Pr}_{Y}(F)\subset{\Pr}_{Y}(\operatorname*{dom}f).\label{r17}\end{equation}
 Indeed, let $y\in{}^{ic}({\Pr}_{Y}(\operatorname*{dom}f))$. Then
$0\in{}^{ic}({\Pr}_{Y}(\operatorname*{dom}f^{\prime}))$ with $f^{\prime}:=f_{(0,y,0,0)}$
because $\operatorname*{dom}f^{\prime}=\operatorname*{dom}f-(0,y,0,0)$.
Since $f^{\prime}\in\mathcal{G}_{s}$, by (i) we get $G(M_{f^{\prime}})=\{(x,x^{\ast})\mid\exists y^{\ast}\ :\
(x,y,x^{\ast},y^{\ast})\in M_{f}\}$ is maximal monotone; in particular $G(M_{f^{\prime}})$ is nonempty,
and so $y\in\Pr_{Y}(F)$. Hence the first inclusion of (\ref{r17})
holds while the second one is obvious.

Because $f\in\mathcal{G}_{s}$, from (\ref{ineq-s}), we have that
$\varphi_{F}\leq f\leq\operatorname*{conv}c_{F}$. It follows that
\[
F\subset\operatorname*{conv}F\subset\operatorname*{dom}(\operatorname*{conv}c_{F})\subset\operatorname*{dom}f\subset\operatorname*{dom}\varphi_{F},\]
 whence \begin{align}
{\Pr}_{Y}(F) & \subset{\Pr}_{Y}(\operatorname*{conv}F)=\operatorname*{conv}({\Pr}_{Y}(F))\subset{\Pr}_{Y}(\operatorname*{dom}(\operatorname*{conv}c_{F}))\notag\\
 & \subset{\Pr}_{Y}(\operatorname*{dom}f)\subset{\Pr}_{Y}(\operatorname*{dom}\varphi_{F}).\label{r12}\end{align}
 This together with Lemma \ref{lem4}~(ii) yield \begin{align}
\operatorname*{aff}({\Pr}_{Y}(F)) & =\operatorname*{aff}({\Pr}_{Y}(\operatorname*{conv}F))\subset\operatorname*{aff}({\Pr}_{Y}(\operatorname*{dom}(\operatorname*{conv}c_{F}))\nonumber \\
 & \subset\operatorname*{aff}({\Pr}_{Y}(\operatorname*{dom}f))\subset\operatorname*{aff}({\Pr}_{Y}(\operatorname*{dom}\varphi_{F}))\subset\overline{\operatorname*{aff}}\,({\Pr}_{Y}(F)).\label{r13}\end{align}

If $\operatorname*{aff}({\Pr}_{Y}(\operatorname*{conv}F))$ $(=\operatorname*{aff}({\Pr}_{Y}(F)))$
is closed, all inclusions in (\ref{r13}) become equalities; hence
\begin{equation}
{}^{{ic}}({\Pr}_{Y}(F))\subset{}^{ic}({\Pr}_{Y}(\operatorname*{conv}F))\subset{}^{ic}({\Pr}_{Y}(\operatorname*{dom}f))\subset{}^{ic}({\Pr}_{Y}(\operatorname*{dom}\varphi_{F}))\label{r15}\end{equation}
 in this case.

Assume that $^{ic}\left(\Pr_{Y}(\operatorname*{dom}f)\right)\neq\emptyset$.
Taking into account (\ref{suppl}) and (\ref{r17}), we know that
$\operatorname*{aff}({\Pr}_{Y}(\operatorname*{dom}f))=\operatorname*{aff}({\Pr}_{Y}(F))$
is closed and ${}^{{ic}}({\Pr}_{Y}(F))={}^{ic}({\Pr}_{Y}(\operatorname*{dom}f))$.
Hence from (\ref{r15}) we obtain \begin{equation}
^{{ic}}({\Pr}_{Y}(F))={}^{ic}({\Pr}_{Y}(\operatorname*{conv}F))={}^{ic}({\Pr}_{Y}(\operatorname*{dom}f)).\label{r18}\end{equation}
 If $^{{ic}}({\Pr}_{Y}(F))\neq\emptyset$ then $\operatorname*{aff}({\Pr}_{Y}(F))$
is closed, and so (\ref{r15}) holds. Hence $^{ic}({\Pr}_{Y}(\operatorname*{dom}f))\neq\emptyset$,
whence (\ref{r18}) holds again. Since $X\times Y\times X^{*}\times Y^{*}$
is a Banach space and $f\in\Gamma_{s}(X\times Y\times X^{\ast}\times Y^{\ast})$,
by \cite[Prop.\ 2.7.2]{Zalinescu:02} we have $^{ic}\left(\Pr_{Y}(\operatorname*{dom}f)\right)=\operatorname*{ri}\left(\Pr_{Y}(\operatorname*{dom}f)\right)$;
taking (\ref{r17}) into account, we get \begin{equation}
^{{ic}}({\Pr}_{Y}(F))={}^{ic}({\Pr}_{Y}(\operatorname*{conv}F))={}^{ic}({\Pr}_{Y}(\operatorname*{dom}f))=\operatorname*{ri}\left({\Pr}_{Y}(\operatorname*{dom}f)\right).\label{r19}\end{equation}
 By Corollary \ref{c-v1} we have that $\varphi_{F}$ is a strong
representative of $F=M_{f}$. Hence, from (\ref{r19}) applied for
$\varphi_{F}$, we find ${}^{ic}({\Pr}_{Y}(\operatorname*{dom}\varphi_{F}))={}^{{ic}}({\Pr}_{Y}(F))$,
thereby completing the proof of (\ref{r14}).\hfill{}$\square$

\medskip{}

For $F:X\times Y\rightrightarrows X^{\ast}\times Y^{\ast}$ and $A:X\rightarrow Y$
a continuous linear operator, we consider $F_{A}:X\times Y\rightrightarrows X^{\ast}\times Y^{\ast}$
defined by \[
\operatorname*{gph}F_{A}:=\{(x,y,x^{\ast},y^{\ast})\in X\times Y\times X^{\ast}\times Y^{\ast}\mid(x^{\ast}-A^{\top}y^{\ast},y^{\ast})\in F(x,Ax+y)\},\]
 where $A^{\top}:Y^{*}\to X^{*}$ is the adjoint of $A$, or $F_{A}(x,y)=B^{\top}FB(x,y)$
with $B(x,y):=(x,y+Ax)$ for $(x,y)\in X\times Y$.

Since $B:X\times Y\rightarrow X\times Y$ is an isomorphism of normed
vector spaces (with $B^{\top}(x^{\ast},y^{\ast})=(x^{\ast}+A^{\top}y^{\ast},y^{\ast})$),
if $F$ is strongly--representable, (maximal) monotone then $F_{A}$
is strongly--representable, (maximal) monotone. Moreover, if $f$
is a (strong) representative of $F$ then $f_{A}:=f\circ L$ is a
(strong) representative of $F_{A}$, where $L:=B\times(B^{-1})^{\top}$.
In an extended form\[
f_{A}(x,y,x^{*},y^{*})=f(x,y+Ax,x^{*}-A^{\top}y^{*},y^{*}),\quad(x,y,x^{*},y^{*})\in X\times Y\times X^{\ast}\times Y^{\ast}.\]
 Note that $y\in{\Pr}_{Y}(\operatorname*{dom}f_{A})$ iff $y-Ax\in{\Pr}_{Y}(\operatorname*{dom}f)$
for some $(x,y)\in\Pr\nolimits _{X\times Y}(\operatorname*{dom}f)$,
$(x,y,x^{*},y^{*})\in M_{f_{A}}$ iff $(x,y+Ax,x^{*}-A^{\top}y^{*},y^{*})\in M_{f}$,
and $(M_{f})_{A}=M_{f_{A}}$, for every $f\in{\cal F}$.

Using the previous result for $F_{A}$ we get the next two consequences.

\begin{corollary} \label{c14}Assume that $X,Y$ are Banach spaces,
$f\in\mathcal{G}_{s}\left(X\times Y\times X^{\ast}\times Y^{\ast}\right)$
and $A\in L(X,Y)$. Then \begin{align*}
{}^{ic}\{y-Ax\mid(x,y)\in\operatorname*{dom}M_{f}\} & ={}^{ic}\{y-Ax\mid(x,y)\in\operatorname*{conv}(\operatorname*{dom}M_{f})\}\\
 & ={}^{ic}\{y-Ax\mid(x,y)\in\Pr\nolimits _{X\times Y}(\operatorname*{dom}f)\}\\
 & =\operatorname*{ri}(\{y-Ax\mid(x,y)\in\operatorname*{dom}M_{f}\}).\end{align*}

Assume that $0\in{}^{ic}\{y-Ax\mid(x,y)\in\Pr\nolimits _{X\times Y}(\operatorname*{dom}f)\}$
(or equivalently $0\in{}^{ic}\{y-Ax\mid(x,y)\in\operatorname*{dom}M_{f}\}$).
Then the multifunction $G(F_{A})$ whose graph is $\{(x,x^{\ast})\in X\times X^{\ast}\mid\exists y^{\ast}\in Y^{\ast}:(x^{\ast}-A^{\top}y^{\ast},y^{\ast})\in M_{f}(x,Ax)\}$
is strongly-representable, a strong representative is given by $\overline{g}$
where $g:X\times X^{\ast}\rightarrow\overline{\mathbb{R}}$ is defined
by\[
g(x,x^{\ast})=\inf\{f(x,Ax,x^{\ast}-A^{\top}y^{\ast},y^{\ast})\mid y^{\ast}\in Y^{\ast}\}\quad\forall(x,x^{\ast})\in X\times X^{\ast};\]
 moreover $G(F_{A})$ is maximal monotone. In fact $G(F_{A})=M_{g}=M_{\overline{g}}=M_{g^{\square}}$
and \[
g^{\square}(x,x^{\ast})=\min\{f^{\square}(x,Ax,x^{\ast}-A^{\top}y^{\ast},y^{\ast})\mid y^{\ast}\in Y^{\ast}\}\quad\forall(x,x^{\ast})\in X\times X^{\ast}.\]

\end{corollary}

\begin{theorem} \label{c15}Assume that $X,Y$ are Banach spaces,
$f\in\mathcal{G}_{s}(X\times X^{\ast})$, $g\in\mathcal{G}_{s}(Y\times Y^{\ast})$
and $A\in L(X,Y)$. Then\begin{align*}
{}^{ic}(\operatorname*{dom}M_{g}-A(\operatorname*{dom}M_{f})) & ={}^{ic}(\operatorname*{conv}(\operatorname*{dom}M_{g}-A(\operatorname*{dom}M_{f})))\\
 & ={}^{ic}\!\left({\Pr}_{Y}(\operatorname*{dom}g)-A({\Pr}_{X}(\operatorname*{dom}f))\right)\\
 & =\operatorname*{ri}\left(\operatorname*{dom}M_{g}-A(\operatorname*{dom}M_{f})\right).\end{align*}

If, in addition, $0\in{}^{ic}(\operatorname*{dom}g-A(\operatorname*{dom}f))$
(or equivalently $0\in{}^{ic}(\operatorname*{dom}M_{g}-A(\operatorname*{dom}M_{f}))$)
then $M_{f}+A^{\top}M_{g}A$ is strongly representable (and maximal
monotone) having as strong representative the function $\overline{k}$,
where\begin{equation}
k:X\times X^{\ast}\rightarrow\overline{\mathbb{R}},\quad\
k(x,x^{\ast}):=\inf\{f(x,x^{\ast}-A^{\top}y^{\ast})+g(Ax,y^{\ast})\mid y^{\ast}\in Y^{\ast}\}.\label{r16}\end{equation}
 Moreover, $M_{f}+A^{\top}M_{g}A=M_{k}=M_{\overline{k}}=M_{k^{\square}}$
and \[
k^{\square}(x,x^{\ast}):=\min\{f^{\square}(x,x^{\ast}-A^{\top}y^{\ast})+g^{\square}(Ax,y^{\ast})\mid y^{\ast}\in Y^{\ast}\}\quad\forall(x,x^{\ast})\in X\times X^{\ast}.\]

\end{theorem}

Proof. Consider $\phi:X\times Y\times X^{\ast}\times Y^{\ast}$ defined
by $\phi(x,y,x^{\ast},y^{\ast}):=f(x,x^{\ast})+g(y,y^{\ast})$. Then
$\phi^{\ast}(x^{\ast},y^{\ast},x^{\ast\ast},y^{\ast\ast})=f^{\ast}(x^{\ast},x^{\ast\ast})+g^{\ast}(y^{\ast},y^{\ast\ast})$,
and so $\phi\in\mathcal{G}_{s}(X\times Y\times X^{\ast}\times Y^{\ast})$.
Moreover, for $F:=M_{\phi}$ we have $G(F_{A})=M_{f}+A^{\top}M_{g}A$.
The conclusion follows using the preceding corollary. \hfill{}$\square$

\medskip{}

Taking $X=Y$ and $A=\operatorname*{Id}_{X}$ in the preceding theorem
we get the following result which shows that the Rockafellar Conjecture
on the sum of maximal monotone operators is true in the strongly-representable
case.

\begin{corollary} \label{c16}Let $X$ be a Banach space and let
$M,N:X\rightrightarrows X^{\ast}$ be strongly representable. Then
$^{ic}(\operatorname*{dom}M-\operatorname*{dom}N)={}^{ic}(\operatorname*{conv}(\operatorname*{dom}M)-\operatorname*{conv}(\operatorname*{dom}N))$
(is a convex set). If $0\in{}^{ic}(\operatorname*{dom}M-\operatorname*{dom}N)$
then $M+N$ is strongly representable; in particular $M+N$ is maximal
monotone. Moreover, $\operatorname*{cl}(\operatorname*{dom}(M+N))$
and $\operatorname*{cl}(\operatorname{Im}(M+N))$ are convex sets.
\end{corollary}

\begin{remark} \label{rem 5}Since every subdifferential is strongly-representable,
the previous corollary together with \cite[Theorem 26.1]{Simons:98}
show that every strongly-representable operator is maximal monotone
locally. \end{remark}

\begin{theorem} \label{t-last} If $X$ is a Banach space, $M:X\rightrightarrows X^{\ast}$
is strongly representable, and $N:X\rightrightarrows X^{\ast}$ is
maximal monotone with $\operatorname*{dom}N=X$, then $M+N$ is maximal
monotone. \end{theorem}

Proof. In order to prove that $M+N$ is maximal monotone we wish to
apply \cite[Th. 3.4]{Voisei-tscr}, that is, to show that $M+N$ is
representable and $\varphi_{M+N}\geq c.$ Since $M+N$ is representable
by \cite[Cor. 5.6]{Voisei-tscr}, we have only to prove that $\varphi_{M+N}\geq c,$
or equivalently that $\overline{x}\in\operatorname*{dom}(M+N)$ whenever
$\overline{z}:=(\overline{x},\overline{x}^{\ast})$ is monotonically
related to $M+N$ (because always for a monotone operator $S:X\rightrightarrows X^{\ast}$
one has $(\operatorname*{dom}S)\times X^{\ast}\subset\lbrack\varphi_{S}\geq c];$
see \cite[Corollary 5.6]{Voisei-tscr}).

According to Corollary \ref{c-v1}, we may choose $f$ to be a strong
representative for $M$ such that $f\in\mathcal{G}_{s\times w^{\ast}}(X\times X^{\ast})$.
Let $\overline{z}=(\overline{x},\overline{x}^{\ast})$ be monotonically
related to $M+N$. Taking $M_{0}:=M-\overline{z}$ and $N_{0}=N-(\overline{x},0),$
then $\operatorname*{gph}(M+N)-\overline{z}=\operatorname*{gph}(M_{0}+N_{0})$
and $(0,0)$ is monotonically related to $M_{0}+N_{0};$ moreover,
$f_{\overline{z}}\in\mathcal{G}_{s\times w^{\ast}}(X\times X^{\ast}),$
$f_{\overline{z}}$ is a strong representative of $M_{0}$ and $\operatorname*{dom}N_{0}=X.$
If we prove that $0\in\operatorname*{dom}(M_{0}+N_{0})$ then $\overline{x}\in\operatorname*{dom}(M+N).$
Hence without loss of generality we assume that $\overline{z}=0,$
and so \begin{equation}
c(u,u^{\ast}+v^{\ast})\geq0\quad\ \forall(u,u^{\ast})\in M,\
(u,v^{\ast})\in N.\label{m.r.t.}\end{equation}

Fix $(x_{0},x_{0}^{\ast})\in\operatorname*{dom}f$ and let $[-x_{0},x_{0}]:=\{tx_{0}\mid-1\leq t\leq1\}$
and $C_{\varepsilon}:=[-x_{0},x_{0}]+\varepsilon U$ for $\varepsilon>0,$
where $U:=U_{X}:=\{x\in X\mid\left\Vert x\right\Vert \leq1\}$. Since
$N$ is locally bounded and $[-x_{0},x_{0}]$ is compact there is
$\varepsilon_{0}>0$ such that $N$ is bounded on $C_{\varepsilon_{0}}$,
that is, there is $K>0$ such that \begin{equation}
\Vert v^{\ast}\Vert\leq K\quad\forall v\in C_{\varepsilon_{0}},\ \forall v^{\ast}\in N(v).\label{lb}\end{equation}
 Take $C:=C_{\varepsilon_{0}/2}$. Notice that $C$ is symmetric,
that is, $-x\in C$ for every $x\in C$.

Let $\psi_{n}(x)=\iota_{C}(x)+\frac{n}{2}\Vert x\Vert^{2}$ for $x\in X$
and $\Psi_{n}(x,x^{\ast})=\psi_{n}(x)+\psi_{n}^{\ast}(x^{\ast})$
for $(x,x^{\ast})\in X\times X^{\ast}$ and $n\geq1$. Since $C$
is symmetric we have $\psi_{n}(x)=\psi_{n}(-x)$ and $\psi_{n}^{\ast}(x^{\ast})=\psi_{n}^{\ast}(-x^{\ast})$
for all $x\in X$ and $x^{\ast}\in X^{\ast};$ it follows that $\Psi_{n}\geq\pm c.$
Moreover, \begin{equation}
\psi_{n}^{\ast}(x^{\ast})=\min\big\{\sigma_{C}(u^{\ast})+\tfrac{1}{2n}\left\Vert x^{\ast}-u^{\ast}\right\Vert ^{2}\mid u^{\ast}\in X^{\ast}\big\}\geq0\quad\forall x^{\ast}\in X^{\ast},\label{psin*}\end{equation}
 and $\psi_{n}^{\ast}$ is finite and continuous on $X^{\ast},$ where
for $A\subset X$ and $x^{\ast}\in X^{\ast},$ $\sigma_{A}(x^{\ast}):=\iota_{A}^{\ast}(x^{\ast})=\sup_{x\in A}\left\langle x,x^{\ast}\right\rangle .$

Since $\Psi_{n}$ is continuous at $(x_{0},x_{0}^{\ast}),$ $f\geq c,$
$\Psi_{n}\geq-c,$ $f^{\ast}\geq c$ and $\Psi_{n}^{\ast}\geq-c,$
as in the proof of Proposition \ref{t-v}, applying the fundamental
duality formula, we get \[
\inf_{w\in X\times X^{\ast}}(f(w)+\Psi_{n}(w))=0\quad\ \forall n\geq1.\]

Therefore, for every $n\geq1$ there is $z_{n}=(x_{n},x_{n}^{\ast})$
such that $f(z_{n})+\Psi_{n}(z_{n})<n^{-2}$. Since $x_{n}\in C$,
we know that $\Vert x_{n}\Vert\leq\Vert x_{0}\Vert+\varepsilon_{0}/2$
for $n\geq1$.

As seen above, $f\geq c$ and $\Psi_{n}\geq-c$; it follows that \begin{equation}
\Psi_{n}(z_{n})+c(z_{n})\leq n^{-2},\quad\
f(z_{n})<c(z_{n})+n^{-2}\quad\ \forall n\geq1.\label{Brr}\end{equation}

From (\ref{Brr}), Corollary \ref{BR} provides $w_{n}=(y_{n},y_{n}^{\ast})\in M$
such that $\Vert w_{n}-z_{n}\Vert<2/n$ for $n\geq1$.

Pick $v_{n}^{\ast}\in N(y_{n})$. Then for $n\geq4/\varepsilon_{0}$
we have that $y_{n}\in C_{\varepsilon_{0}},$ and so $\Vert v_{n}^{\ast}\Vert\leq K$.
Using (\ref{z18}), this yields\begin{align*}
\frac{n}{2}\Vert x_{n}\Vert^{2}+\psi_{n}^{\ast}(x_{n}^{\ast}) & =\Psi_{n}(z_{n})\leq-c(z_{n})+n^{-2}\leq-c(w_{n})+\left\vert c(w_{n})-c(z_{n})\right\vert +n^{-2}\\
 & \leq-c(w_{n})+\tfrac{1}{2}\Vert w_{n}-z_{n}\Vert^{2}+\Vert z_{n}\Vert\Vert z_{n}-w_{n}\Vert+n^{-2}\\
 & \leq-c(w_{n})+2n^{-1}\Vert z_{n}\Vert+3n^{-2}\leq-c(w_{n})+2n^{-1}\Vert x_{n}\Vert+2n^{-1}\Vert x_{n}^{\ast}\Vert+3n^{-2}\\
 & =-c(y_{n},y_{n}^{\ast}+v_{n}^{\ast})+c(y_{n},v_{n}^{\ast})+2n^{-1}\Vert x_{n}\Vert+2n^{-1}\Vert x_{n}^{\ast}\Vert+3n^{-2}\\
 & \leq K(\Vert x_{n}\Vert+2n^{-1})+2n^{-1}\Vert x_{n}\Vert+2n^{-1}\Vert x_{n}^{\ast}\Vert+3n^{-2}\\
 & \leq K\Vert x_{n}\Vert+2n^{-1}\Vert x_{n}^{\ast}\Vert+Ln^{-1}\end{align*}
 for $n\geq4/\varepsilon_{0}$, where $L:=2K+2\Vert x_{0}\Vert+\varepsilon_{0}+3$.
Hence, for $n\geq4/\varepsilon_{0}$ we have\[
\frac{n}{2}\Vert x_{n}\Vert^{2}-K\Vert x_{n}\Vert+[\psi_{n}^{\ast}(x_{n}^{\ast})-2n^{-1}\Vert x_{n}^{\ast}\Vert-Ln^{-1}]\leq0,\]
 or equivalently\begin{equation}
\tfrac{1}{2}(\Vert nx_{n}\Vert-K)^{2}+[n\psi_{n}^{\ast}(x_{n}^{\ast})-2\Vert x_{n}^{\ast}\Vert]\leq\tfrac{1}{2}K^{2}+L.\label{estimare}\end{equation}

We claim that\begin{equation}
n\psi_{n}^{\ast}(x^{\ast})\geq3\Vert x^{\ast}\Vert-18\quad\ \forall x^{\ast}\in X^{\ast},~\forall n\geq6/\varepsilon_{0}.\label{psi*}\end{equation}

The condition $n\geq6/\varepsilon_{0}$ implies $nC\supset3U$; whence
$n\sigma_{C}(u^{\ast})=\sigma_{nC}(u^{\ast})\geq3\Vert u^{\ast}\Vert$
for every $u^{\ast}\in X^{\ast}$.

For fixed $x^{\ast}\in X^{\ast}$ we consider two cases: a) $\Vert x^{\ast}-u^{\ast}\Vert<6$
and b) $\Vert x^{\ast}-u^{\ast}\Vert\geq6$.

\noindent If a) holds then $\Vert u^{\ast}\Vert\geq\Vert x^{\ast}\Vert-\Vert x^{\ast}-u^{\ast}\Vert>\Vert x^{\ast}\Vert-6$
and so\[
\sigma_{nC}(u^{\ast})+\tfrac{1}{2}\Vert x^{\ast}-u^{\ast}\Vert^{2}\geq\sigma_{nC}(u^{\ast})\geq3\Vert u^{\ast}\Vert\geq3\Vert x^{\ast}\Vert-18.\]

\noindent If b) holds then $\frac{1}{2}\Vert x^{\ast}-u^{\ast}\Vert^{2}\geq3\Vert x^{\ast}-u^{\ast}\Vert$
and so\[
\sigma_{nC}(u^{\ast})+\tfrac{1}{2}\Vert x^{\ast}-u^{\ast}\Vert^{2}\geq3\Vert u^{\ast}\Vert+3\Vert x^{\ast}-u^{\ast}\Vert\geq3\Vert x^{\ast}\Vert.\]
 Taking into account (\ref{psin*}) we obtain that our claim is true.
Using (\ref{psi*}), from (\ref{estimare}) we get \[
\tfrac{1}{2}(\Vert nx_{n}\Vert-K)^{2}+\left\Vert x_{n}^{\ast}\right\Vert \leq\tfrac{1}{2}K^{2}+L+18\quad\forall n\geq6/\varepsilon_{0}.\]
 Hence necessarily $\left\Vert x_{n}\right\Vert \rightarrow0$ and
$(x_{n}^{\ast})$ is bounded. On a subnet, denoted for simplicity
by the same index, $x_{n}^{\ast}\rightarrow x^{\ast}$ weakly-star
in $X^{\ast}$. Passing to limit in (\ref{Brr}) we get $(0,x^{\ast})\in\lbrack f=c]=M$
and so $\overline{x}=0\in\operatorname*{dom}M=\operatorname*{dom}(M+N)$.
The proof is complete. \hfill{}$\square$

\strut

The previous theorem allows us to recover the results in \cite[Theorem 42.2]{Simons:98}
and its extension \cite[Corollary 2.9(a)]{Veronas-2000 reg max mon and sum}.

\begin{corollary}If $X$ is a Banach space, $\varphi\in\Gamma_{s}(X)$
and $L:X\rightarrow X^{*}$ is linear positive then $\partial\varphi+L$
is maximal monotone.

\end{corollary}

\begin{corollary}If $X$ is a Banach space, $\varphi\in\Gamma_{s}(X)$
and $N:X\rightrightarrows X^{*}$ is maximal monotone with $\operatorname*{dom}N=X$
then $\partial\varphi+N$ is maximal monotone.

\end{corollary}


\section{Comparison with other classes of operators}

Recall that $M\subset Z:=X\times X^{\ast}$ is called \emph{locally
maximal monotone} if for every open convex $U$ in $X^{*}$ such that
$U\cap\operatorname*{Im}M\neq\emptyset$ and $z\in X\times U\setminus\operatorname*{gph}M$
there is $w\in\operatorname*{gph}M\cap X\times U$ such that $c(z-w)<0$
(see \cite{Fitz-phelps - 92 ba}).

\begin{theorem}Every strongly-representable operator is locally maximal
monotone.

\end{theorem}

Proof. Let $M$ be a strongly-representable operator with a strong-representative
$f\in{\cal G}_{s\times w^{*}}(X\times X^{*})$. According to \cite[Proposition 3.2]{Fitz-phelps - 92 ba},
it suffices to prove the counter-positive form of the definition on
bounded convex weakly star closed sets, that is, for every weakly-star
closed bounded convex set $C$ in $X^{*}$ such that $C\cap\operatorname*{Im}M\neq\emptyset$
and $z\in X\times\operatorname*{int}C$ with $c(z-w)\ge0$, for all
$w\in\operatorname*{gph}M\cap X\times C$ then $z\in\operatorname*{gph}M$.

Without loss of generality we may assume that $z=0$, whence $0\in\operatorname*{int}C$;
therefore, there is $r_{0}>0$ such that $2r_{0}U\subset C$, where
$U:=\{x^{*}\in X^{*}\mid\|x^{*}\|\le1\}$. Thus \begin{equation}
c(w)\ge0\quad\forall w\in\operatorname*{gph}M\cap(X\times C).\label{r-m}\end{equation}
 Let $\psi_{r}(x,x^{*})=r\|x\|+\iota_{rU}(x^{*})$ for $(x,x^{*})\in X\times X^{*}$
and $0<r\le r_{0}$.

Let us fix $r\in(0,r_{0})$. As previously seen, from the fundamental
duality formula and from $f\ge c$, $f^{*}\ge c$, $\psi_{r}\ge\pm c$,
$\psi_{r}^{*}\ge\pm c$ we get $\inf(f+\psi_{r})=0$, and this implies
the existence of $z_{n}=(x_{n},x_{n}^{*})$ such that $f(z_{n})+\psi_{r}(z_{n})<n^{-2}$
for $n\ge1$. Again, because $f\geq c$ and $\psi_{r}\geq-c$, we
get $\psi_{r}(z_{n})+c(z_{n})\le n^{-2}$ and $f(z_{n})<c(z_{n})+n^{-2}$
for $n\ge1$. Corollary \ref{BR} provides $w_{n}\in M$ such that
$\|w_{n}-z_{n}\|<2/n$ for $n\ge1$. Note that $w_{n}\in\operatorname*{gph}M\cap(X\times C)$,
and so $c(w_{n})\ge0$ for every $n\ge2/r$. Taking into account (\ref{r-m})
and (\ref{z18}) we get \begin{align*}
r\|x_{n}\| & =\psi_{r}(z_{n})\le-c(z_{n})+n^{-2}\le-c(w_{n})+|c(z_{n})-c(w_{n})|+n^{-2}\\
 & \leq\tfrac{1}{2}\left\Vert w_{n}-z_{n}\right\Vert ^{2}+\left\Vert z_{n}\right\Vert \cdot\left\Vert w_{n}-z_{n}\right\Vert +n^{-2}\leq2n^{-1}\left\Vert x_{n}\right\Vert +2rn^{-1}+3n^{-2},\end{align*}
 for $n\ge2/r$. This inequality shows that $x_{n}\rightarrow0$ strongly
in $X$, as $n\rightarrow\infty$. Since $\|x_{n}^{*}\|\le r$ for
$n\ge1$, on a subnet, denoted for simplicity by the same index, $x_{n}^{*}\rightarrow\overline{x}_{r}^{*}$
weakly-star in $X^{*}$; hence $\|\overline{x}_{r}^{*}\|\le r$ and
$z_{n}\rightarrow(0,\overline{x}_{r}^{*})$ for the topology $s\times w^{*}$
in $X\times X^{*}$. Passing to limit in $c(z_{n})\le f(z_{n})<c(z_{n})+n^{-2}$
we find that $(0,\overline{x}_{r}^{*})\in[f=c]=M$.

We proved that for every $0<r\le r_{0}$, there is $\overline{x}_{r}^{*}$
with $\|\overline{x}_{r}^{*}\|\le r$ such that $\overline{x}_{r}^{*}\in M(0)$.
Since $\overline{x}_{r}^{*}\rightarrow0$, strongly in $X^{*}$, as
$r\searrow0$ and $M$ is maximal monotone thus it has closed values,
it yields that $0\in M(0)$, i.e., $z=0\in\operatorname*{gph}M$.
The proof is complete. \hfill{}$\square$

\strut

Using a different argument the previous result allows us to recover
the convexity of the closure of the range of a strongly-representable
operator (see \cite[Theorem 3.5]{Fitz-phelps - 92 ba}).

\strut

Recall that $M\subset Z:=X\times X^{\ast}$ is called NI if\[
\inf_{(u,u^{\ast})\in M}\langle u^{\ast}-x^{\ast},\widehat{u}-x^{\ast\ast}\rangle\leq0\quad\
\forall(x^{\ast},x^{\ast\ast})\in Z^{\ast},\]
 or equivalently $\Phi_{M}(x^{\ast},x^{\ast\ast}):=c_{M}^{*}(x^{\ast},x^{\ast\ast})\geq\langle x^{\ast},x^{\ast\ast}\rangle$
for every $(x^{\ast},x^{\ast\ast})\in X^{\ast}\times X^{\ast\ast}$.

\begin{proposition} \label{NI} Let $M\subset X\times X^{\ast}$
be maximal monotone and NI. Then $M$ is strongly-representable. \end{proposition}

Proof. Since $M$ is maximal monotone we have that $c_{M}\ge\varphi_{M}\ge c$
and $M=M_{\varphi_{M}}$. It follows that $\varphi_{M}^{*}\ge c_{M}^{*}$.
Since $M$ is NI we have that $c_{M}^{*}\ge c$, and so $\varphi_{M}\in\mathcal{G}_{s}(Z)$.
Hence $\varphi_{M}$ is a strong representative of $M$. \hfill{}$\square$

\medskip{}

For skew bounded operators the converse of Proposition \ref{NI} holds.

\begin{corollary} Let $S:X\rightrightarrows X^{*}$ be skew, that
is, $\operatorname*{gph}S$ is a linear subspace and $\langle x,x^{*}\rangle=0$
for all $x\in\operatorname*{dom}S$ and $x^{*}\in S(x)$. Consider
the conditions:

(i) $S$ is maximal monotone and NI,

(ii) $S$ is $s\times w^{*}$--closed in $X\times X^{*}$ and $S^{*}$
is monotone in $X^{**}\times X^{*}$,

(iii) $S$ is strongly-representable.

Then (i) $\Leftrightarrow$ (ii) $\Rightarrow$ (iii). If in addition
$S:X\rightarrow X^{*}$ has $\operatorname*{dom}S=X$, then (iii)
$\Rightarrow$ (ii). Here $(x^{**},x^{*})\in\operatorname*{gph}S^{*}$
iff $\langle u,x^{*}\rangle=\langle u^{*},x^{**}\rangle$ for every
$(u,u^{*})\in\operatorname*{gph}S$. \end{corollary}

Proof. The implication (i) $\Rightarrow$ (ii) follows from the proof
of Proposition \ref{NI} and the fact that $S$ is skew. In this case
$\psi_{S}:=\varphi_{S}^{\square}=\iota_{S}$ is a strong-representative
of $S$ with $\iota_{S}^{*}=\iota_{(-S^{*})^{-1}}$. The implication
(ii) $\Rightarrow$ (i) follows from $\Phi_{S}=\iota_{S}^{*}=\iota_{(-S^{*})^{-1}}$.
For (ii) $\Rightarrow$ (iii) notice that $\iota_{S}$ is a strong-representative
of $S$.

(iii) $\Rightarrow$ (ii) Let $f\in\mathcal{G}_{s}$ be a strong-representative
of $S$. Since $S$ is maximal monotone and skew $f\ge\varphi_{S}=\iota_{S}$;
hence $\iota_{S}^{*}=\iota_{(-S^{*})^{-1}}\ge f^{*}\ge c$, that is,
$S^{*}$ is monotone in $X^{**}\times X^{*}$ and $S$ is $s\times w^{*}-$closed
in $X\times X^{*}$, since $\operatorname*{cl}_{s\times w^{*}}S$
remains skew. \hfill{}$\square$

\begin{corollary} \label{c-z4}Let $f\in\mathcal{G}_{s}(Z)$. For
every $(x,x^{\ast})\in X\times X^{\ast}$ and every $\varepsilon>0$
there exists $(x_{\varepsilon},x_{\varepsilon}^{\ast})\in M_{f}$
such that $\left\{ (x_{\varepsilon},x_{\varepsilon}^{\ast})\mid\varepsilon>0\right\} $
is bounded and\[
\left\Vert x-x_{\varepsilon}\right\Vert ^{2}+2\left\langle x-x_{\varepsilon},x^{\ast}-x_{\varepsilon}^{\ast}\right\rangle +\left\Vert x^{\ast}-x_{\varepsilon}^{\ast}\right\Vert ^{2}\leq\varepsilon.\]

\end{corollary}

Proof. Replacing if necessary $f$ by $f_{(x,x^{\ast})}$, we may
(and we do) assume that $(x,x^{\ast})=(0,0)$. As seen in the proof
of Theorem \ref{SR-MM}, $f+h$ is (strongly) coercive. Hence there
exists $r>0$ such that $\{z\in Z\mid f(z)+h(z)\leq1\}\subset rU_{Z}$.
For $\varepsilon\in(0,1]$ take $\varepsilon^{\prime}\in(0,\varepsilon)$
such that $11\varepsilon^{\prime}+6r\sqrt{2\varepsilon^{\prime}}=\varepsilon$.
Since $\inf(f+h)=0$, there exists $w_{\varepsilon}\in Z$ such that
$f(w_{\varepsilon})+h(w_{\varepsilon})\leq\varepsilon^{\prime}$.
Since $f-c\geq0$ and $h\geq-c$, it follows that \[
f(w_{\varepsilon})\leq c(w_{\varepsilon})+\varepsilon^{\prime},\quad\tfrac{1}{2}\left\Vert w_{\varepsilon}\right\Vert ^{2}+c(w_{\varepsilon})\leq\varepsilon^{\prime}.\]
 Using now Theorem \ref{t-z5}, for $\varepsilon^{\prime}>0$ and
$w_{\varepsilon}$ we get $z_{\varepsilon}\in M_{f}$ such that $\Vert w_{\varepsilon}-z_{\varepsilon}\Vert\leq\delta:=\sqrt{4.5\varepsilon^{\prime}}$.
Using (\ref{z18}) we get\begin{align*}
\left\Vert z_{\varepsilon}\right\Vert ^{2} & \leq\left(\left\Vert w_{\varepsilon}\right\Vert +\left\Vert z_{\varepsilon}-w_{\varepsilon}\right\Vert \right)^{2}\leq\left\Vert w_{\varepsilon}\right\Vert ^{2}+2r\left\Vert z_{\varepsilon}-w_{\varepsilon}\right\Vert +\left\Vert z_{\varepsilon}-w_{\varepsilon}\right\Vert ^{2},\\
c(z_{\varepsilon}) & =c(w_{\varepsilon})+\left\langle w_{\varepsilon},z_{\varepsilon}-w_{\varepsilon}\right\rangle +c(z_{\varepsilon}-w_{\varepsilon})\leq c(w_{\varepsilon})+r\left\Vert z_{\varepsilon}-w_{\varepsilon}\right\Vert +\tfrac{1}{2}\left\Vert z_{\varepsilon}-w_{\varepsilon}\right\Vert ^{2}.\end{align*}
 Therefore, \[
\left\Vert z_{\varepsilon}\right\Vert ^{2}+2c(z_{\varepsilon})\leq\left\Vert w_{\varepsilon}\right\Vert ^{2}+2c(w_{\varepsilon})+4r\delta+2\delta^{2}\leq2\varepsilon^{\prime}+4r\delta+2\delta^{2}=11\varepsilon^{\prime}+6r\sqrt{2\varepsilon^{\prime}}=\varepsilon.\]
 Taking $z_{\varepsilon}:=z_{1}$ for $\varepsilon\geq1$, the proof
is complete.

\medskip{}

The next result shows that every strongly-representable operator is
of type ANA (see \cite{Simons:98} for this notion).

\begin{corollary} \label{c-z5}Let $f\in\mathcal{G}_{s}(Z)$. Then
for every $(x,x^{\ast})\in X\times X^{\ast}\setminus M_{f}$ there
exists a (bounded) sequence $\left((x_{n},x_{n}^{\ast})\right)_{n\geq1}\subset M_{f}$
such that $x_{n}\neq x$, $x_{n}^{\ast}\neq x^{\ast}$ for every $n$
and \[
\lim\frac{\left\langle x_{n}-x,x_{n}^{\ast}-x^{\ast}\right\rangle }{\left\Vert x_{n}-x\right\Vert \cdot\left\Vert x_{n}^{\ast}-x^{\ast}\right\Vert }=-1.\]

\end{corollary}

Proof. Let $(x,x^{\ast})\in X\times X^{\ast}\setminus M_{f}$. Fix
$(\varepsilon_{n})\subset(0,\infty)$ with $\varepsilon_{n}\rightarrow0$.
Using Corollary \ref{c-z4} we get a bounded sequence $\left((x_{n},x_{n}^{\ast})\right)_{n\geq1}\subset M_{f}$
such that \begin{equation}
\left\Vert x-x_{n}\right\Vert ^{2}+2\left\langle x-x_{n},x^{\ast}-x_{n}^{\ast}\right\rangle +\left\Vert x^{\ast}-x_{n}^{\ast}\right\Vert ^{2}\leq\varepsilon_{n}^{2}\quad\forall n\in\mathbb{N}.\label{z16}\end{equation}
 Hence \begin{equation}
\left\vert \left\Vert x-x_{n}\right\Vert -\left\Vert x^{\ast}-x_{n}^{\ast}\right\Vert \right\vert \leq\varepsilon_{n}\quad\forall n\in\mathbb{N}.\label{z15}\end{equation}
 Then there exist $\gamma>0$ and $n_{0}\in\mathbb{N}$ such that
$\left\Vert x-x_{n}\right\Vert \geq2\gamma$ for all $n\geq n_{0}$.
Otherwise for some increasing sequence $(n_{k})\subset\mathbb{N}$
we have $x_{n_{k}}\rightarrow x$. Using (\ref{z15}) we get $x_{n_{k}}^{\ast}\rightarrow x^{\ast}$.
This yields the contradiction\[
\left\langle x,x^{\ast}\right\rangle <f(x,x^{\ast})\leq\liminf f(x_{n_{k}},x_{n_{k}}^{\ast})=\lim\left\langle x_{n_{k}},x_{n_{k}}^{\ast}\right\rangle =\left\langle x,x^{\ast}\right\rangle .\]
 From (\ref{z15}) we obtain \[
\left\vert \frac{\left\Vert x^{\ast}-x_{n}^{\ast}\right\Vert }{\left\Vert x-x_{n}\right\Vert }-1\right\vert \leq\frac{\varepsilon_{n}}{2\gamma}\quad\forall n\geq n_{0},\]
 whence $\lim\left\Vert x^{\ast}-x_{n}^{\ast}\right\Vert /\left\Vert x-x_{n}\right\Vert =1$.
Hence $\left\Vert x^{\ast}-x_{n}^{\ast}\right\Vert \geq\gamma$ for
$n\geq n_{1}$ for some $n_{1}\geq n_{0}$ and $\lim\left\Vert x-x_{n}\right\Vert /\left\Vert x^{\ast}-x_{n}^{\ast}\right\Vert =1$.
We may assume that $n_{1}=0$. From (\ref{z16}) we get\[
-2\leq\frac{2\left\langle x_{n}-x,x_{n}^{\ast}-x^{\ast}\right\rangle }{\left\Vert x_{n}-x\right\Vert \cdot\left\Vert x_{n}^{\ast}-x^{\ast}\right\Vert }\leq\frac{\varepsilon_{n}^{2}}{\gamma^{2}}-\frac{\left\Vert x-x_{n}\right\Vert }{\left\Vert x^{\ast}-x_{n}^{\ast}\right\Vert }-\frac{\left\Vert x^{\ast}-x_{n}^{\ast}\right\Vert }{\left\Vert x-x_{n}\right\Vert }\quad\forall n\in\mathbb{N},\]
 whence the conclusion follows.

\end{document}